# A New Class of Scaling Matrices for Scaled Trust Region Algorithms

Aydin Ayanzadeh[*,1], Shokoufeh Yazdanian[2], Ehsan Shahamatnia[3]

[1] Department of Applied Informatics, Informatics Institute, Istanbul Technical University, Istanbul, Turkey
[2] Department of Computer science, University of Tabriz, Tabriz, Iran
[3] Department of Computer and Electrical Engineering, Universidade Nova de Lisboa, Lisbon, Portugal

*Abstract*— A new class of affine scaling matrices for the interior point Newton-type methods is considered to solve the nonlinear systems with simple bounds. We review the essential properties of a scaling matrix and consider several well-known scaling matrices proposed in the literature. We define a new scaling matrix that is the convex combination of these matrices. The proposed scaling matrix inherits those interesting properties of the individual matrices and satisfies additional desired requirements. The numerical experiments demonstrate the superiority of the new scaling matrix in solving several important test problems.

*Keywords*— Diagonal Scaling; Bound-Constrained Equations; Scaling Matrices; Trust Region Methods; Affine Scaling.

## I. INTRODUCTION

Diverse applications including signal processing [36] and compressive sensing [3,4] incorporate many convex nonlinear optimization problems. Although specialized algorithms have been developed for some of these problems, the interior-point methods are still the main tool to tackle them. These methods require a feasible initial point [19]. It should be paid attention that proper data collection plays an important role in assessing the results obtained [39]. In addition, a vast variety of soft-computing techniques such as evolutionary computing methods [2, 18, 28, 37, 38] and neural networks [1, 31] include optimization problems, which are sensitive to the initial points. Like verifying solution uniqueness conditions, these tasks convert into linear feasibility problems with strict inequalities or nonlinear feasibility problems with bound constraints [16, 34]. The latter is often challenging so that the existing algorithms need to be theoretically and computationally improved. The nonlinear minimization problem with bound constraints is:

$$F(x)=0, \quad x \in \Omega = \{x \in \mathbb{R}^n \mid l_i \leq x_i \leq u_i; \forall i=1,...,n\}. \quad (1)$$

where $F: X \mapsto \mathbb{R}^n$ is a continuously differentiable mapping, $X \subseteq \mathbb{R}^n$ is an open set containing the n-dimensional box $\Omega = \{x \in \mathbb{R}^n \mid l \leq x \leq u\}$. The vectors $l \in (\mathbb{R} \cup \{-\infty\})^n, u \in (\mathbb{R} \cup \{+\infty\})^n$ are lower and upper bounds on the variables such that $\Omega$ has a nonempty interior.

Efficient methods for the solution of this problem with good local convergence behavior have been proposed. The affine scaling trust region approach forms a practical framework for smooth and nonsmooth box constrained systems of nonlinear equations [6-9]. These kinds of methods use ellipsoidal trust region defined by a diagonal scaling matrix [12]. The diagonal scaling handles the bounds while at each iteration a quadratic model of the object function $\frac{1}{2}\|F\|^2$ is minimized within a trust region around the current iteration.

The main motivation for the current work is a series of papers by Bellavia et al [6-11]. The methods they introduced, STRN and CODOSOL have very good numerical properties [6,7]. These methods are widely used in practice [14, 35] and their efficiency has been proved in several papers [7, 23]. In [8] the authors studied global and fast convergence of an inexact dogleg method. They did not investigate the choice of a suitable scaling matrix and only reported the preliminary results. Later in [11] they focused on medium scale problems and replaced the inexact Newton step by the exact solution of the Newton equation. They considered several diagonal scaling matrices and showed the assumptions required to ensure the convergence. The name of the method is Constrained Dogleg *(CoDo)* method which is freely accessible through the website http://codosol.de.unifi.it. The effectiveness of *(CoDo)* is verified by comparing it to *STRSCNE* [7] and to IATR [9].

---

[*] Corresponding author can be contacted via the journal website.

In these scaling-based algorithms, the performance is influenced by the selection of a scaling matrix. We introduce a new class of scaling matrices, which is obtained by the convex combination of current known matrices. We analyze the numerical performance of CoDoSol (The Matlab Solver of CoDo) for different convex combinations of the scaling matrices and compare the results using performance profile approach. We also use a Projected Affine Scaling Interior Point algorithm to check the local convergence properties of the new scaling matrix.

In section 2 we explain the role of the scaling matrices. In section 3 we consider several scaling matrices and review the requirements and assumptions. In section 4, we introduce a new class of scaling matrices and prove that the new matrices satisfy the required assumptions. Finally, in section 5 we report our conclusion of computational results.

## II. SCALED TRUST REGIONS

In this section we describe the idea behind using scaled matrices to solve problem (1). First note that (1) is closely related to the box constrained optimization problem:

$$\min_{x \in W} f(x) = \frac{1}{2} \| F(x) \|^2 \quad s.t. \ x \in W. \quad (2)$$

Every solution of (1) is a global minimum of (2.2) and if $x^*$ is a minimum of (2) such that $f(x^*) = 0$, then $x^*$ is a solution of (1). The first order optimality conditions of (2) are equivalent to the nonlinear system of equations $D(x) \nabla f(x) = 0, \quad x \in \Omega$, where $\nabla f(x) = F'(x)^T F(x)$ and $D$ is a suitable scaling matrix of order n

$$D(x) = diag(d_1(x),...,d_n(x)). \quad (3)$$

Coleman and Li [12] considered only one choice of the scaling matrix, then Heinkenschloss et al. [21] noted that the optimality conditions holds for a general class of scaling matrices satisfying the conditions:

$$d_i(x) \begin{cases} = 0, & \text{if } x_i = l_i \text{ and } [\nabla f(x)]_i > 0 \\ = 0, & \text{if } x_i = u_i \text{ and } [\nabla f(x)]_i < 0 \\ \geq 0, & \text{if } x_i \in l_i, u_i \text{ and } [\nabla f(x)]_i = 0 \\ > 0, & \text{else} \end{cases} \quad (4)$$

for all $i=1,...,n$ and all $x \in W$. In affine scaling methods in order to handle the bounds the direction of the scaled gradient $\widehat{g_k}$ is defined by $\widehat{g_k} = -D_k \nabla f_k$.

Given an iterate $x_k \in int(W)$ and the trust region size $D_k > 0$, the trust region subproblem for (2) is:

$$\min_{p \in \mathbb{R}^n} m_k(p) \ ; \ \| G_k p \| \leq \Delta_k \ , \ x_k + p \in int(\Omega) \quad (5)$$

where $m_k$ is the norm of the linear model for F(x) at $x_k$, i.e. $m_k(p) = \| F_k + F_k' p \|$ and $G_k = G(x_k) \in \mathbb{R}^{n \times n}$ with $G : \mathbb{R}^n \mapsto \mathbb{R}^{n \times n}$. For $G_k = I$ the standard spherical trust region and for $G_k = D_k^{-\frac{1}{2}}$ the elliptical trust region achieves. In order to solve this subproblem different approaches have been proposed like STRN [6] which combines ideas from the classical trust-region Newton method for unconstrained nonlinear equations and the interior affine scaling approach for constrained optimization problems or CoDoSol [11] which is based on a dogleg procedure tailored for constrained problems.

## III. SCALING MATRICES

We consider the following well known matrices:

- $D^{CL}(x)$ given by Coleman and Li [12]. The diagonal elements are:

$$d_i^{CL}(x) = \begin{cases} u_i - x_i & \text{if } (\nabla f(x))_i < 0 \text{ and } u_i < \infty \\ x_i - l_i & \text{if } (\nabla f(x))_i > 0 \text{ and } l_i > -\infty \\ \min\{x_i - l_i, u_i - x_i\} & \text{if } (\nabla f(x))_i = 0 \\ & \text{and } l_i > -\infty \text{ or } u_i < \infty \\ 1 & \text{Otherwise} \end{cases} \quad (6)$$

- $D^{HUU}(x)$ given by Heinkenschloss et al [21]:

$$d_i^{HUU}(x) = \begin{cases} d_i^{CL}(x) & \text{if } |\nabla f(x)_i| < \min\{x_i - l_i, u_i - x_i\}^p \text{ or} \\ & \min\{x_i - l_i, u_i - x_i\} < |\nabla f(x)_i|^p \\ 1 & \text{Otherwise} \end{cases} \quad (7)$$

where $p > 1$ is a fixed constant.

- $D^{KK}$ given by Kanzow and Klug [26, 18]:

$$d_i^{KK}(x) = \begin{cases} 1 & if\, l_i = -\infty \text{ and } u_i = \infty \\ \min\{x_i - l_i + \\ \quad \gamma \max\{0, -\nabla f(x)_i\}, \\ \quad u_i - x_i + \\ \quad \gamma \max\{0, \nabla f(x)_i\}\} & Otherwise \end{cases} \quad (8)$$

For a given constant $g > 0$. - $D^{HMZ}(x)$ by Hager et al [20]:

$$d_i^{HMZ}(x) = \frac{\chi_i(x)}{\alpha(x)\chi_i(x) + |\nabla f(x)_i|} \quad (9)$$

where

$$\chi_i(x) = \begin{cases} u_i - x_i & if\, \nabla f(x)_i < 0 \text{ and } u_i < \infty \\ x_i - l_i & if\, \nabla f(x)_i > 0 \text{ and } l_i > -\infty \\ 1 & if\, \nabla f(x)_i = 0 \end{cases} \quad (10)$$

and $a(x)$ is a continuous function, strictly positive for any x and uniformly bounded away from zero.

This matrix has been used by authors in study of a cyclic Barzilai-Borwein gradient method for bound constrained minimization problems they replaced the Hessian of the objective function with $/_k I$, where $/_k$ is the classical Barzilai-Borwein parameter. Then, in order to compute the new iterate, they move along the scaled gradient $d_k = -D^{HMZ}(x_k)\nabla f(x_k)$, with $a(x_k) = /_k$.

Bellavia et al in [11] introduced the following assumptions for a scaling matrix and verified that all above scaling matrices satisfy almost all the requirements specified by the following assumption.

A. Assumption
(i) D(x) satisfies (4),
(ii) D(x) is bounded in $W \cap B_r(x)$ for any $x \in \Omega$ and $r > 0$,
(iii) There exists a $\bar{\lambda}$ such that the step size $/_k$ to the boundary from $x_k$ along $\widehat{g_k}$ satisfies $/_k > \bar{/}$ whenever $\|\nabla f_k\|$ is uniformly bounded above,
(iv) For any $\bar{x}$ in $int(W)$ there exist two positive constants $\bar{\rho}$ and $C_{\bar{x}}$ such that $B_r(\bar{x}) \subset int(W)$ and $\|D(x)^{-1}\| £ C_{\bar{x}}$ for any x in $B_{\bar{\rho}/2}(\bar{x})$.

We define the convex combination of the above scaling matrices as a new class of scaling matrices as follows:

$$D^{CON} = a_1 D^{CL} + a_2 D^{HUU} + a_3 D^{KK} + a_4 D^{HMZ}, \quad \begin{matrix} 0 \leq a_i \leq 1 \\ i = 1,2,3,4 \end{matrix}, \quad (11)$$

$$\sum_{i=1}^{4} a_i = 1.$$

This scaling matrix demonstrates the advantages of all the scaling matrices involved in its definition and seems an appropriate matrix with combined properties. we first have to verify that this matrix satisfies four desired requirements as follows

(i) Clearly $D^{CON}$ satisfies (4).

**Table 1.** Test problems.

| Pb # | Name | Dimension | Box |
|---|---|---|---|
| 1 | Bullard-Biegler system [[10], 14.1.3] | 2 | - |
| 2 | Ferraris-Tronconi system [[10], 14.1.4] | 2 | - |
| 3 | Brown's almost linear system [[10], 14.1.5] | 5 | [-2,2] |
| 4 | Robot kinematics problem [[10], 14.1.6] | 8 | [-1,1] |
| 5 | Series of CSTRs, R = .935 [[10], 14.1.8] | 2 | [0,1] |
| 6 | Series of CSTRs, R = .995 [10] 14.1.8] | 2 | [-inf,inf] |
| 7 | Chemical equilibrium system [[21], system 1] | 10 | - |
| 8 | Problem HS34[16] | 3 | [0,100] |
| 9 | Problem Wachter and Biegler [24] | 3 | [0,inf] |
| 10 | Effati-Grosan 1, a = 2 [23] | 2 | [-2,2] |
| 11 | Effati-Grosan 1, a =100 [23] | 2 | [-100,100] |
| 12 | Effati-Grosan 2, a = 2 [23] | 2 | [-2,-2] |
| 13 | Effati-Grosan 2, a =100 [23] | 2 | [-100,100] |
| 14 | Trigexp1 [20] | 1000 | [-100,100] |
| 15 | Troesch [20] | 500 | [-1,1] |

(ii) As Bellavia et al. [11] showed, the above three matrices are bounded in $W \cap B_r(x)$ for $x \in W$ and $r > 0$. This implies that the combination of these matrices is also bounded.

(iii) It has been shown in [11] that all the matrices except $D^{HUU}$ join this property. Thus if we assume $a_2 = 0$ then $D^{CON}$ satisfies this condition.

(iv) Same as before, since all the matrices satisfy condition (iv), the convex combination also verifies this property.

It is worth mentioning that $D^{CL}$ is discontinuous at points where there exists an index $i$ for which $\nabla f(x)_i = 0$ but the matrices $D^{KK}$ and $D^{HUU}$ are locally Lipschitz continuous and continuous respectively. The new matrix $D^{CON}$ is continuous and so we can have a matrix having all the advantages of 3 matrices, a matrix that can be as fast and efficient as the first one and as smooth and continuous as second one.

IV. NUMERICAL EXPERIMENTS

We ran two experiments, CoDoSol on 15 problems and Projected Affine Scaling Interior Point on 2 problems.

### A. Experiment 1

In this section, we apply the CoDoSol algorithm to several problems using 3 scaling matrices $D^{CL}, D^{HUU}, D^{KK}$ and 4 combined matrices:

$$\frac{1}{2}D^{CL} + \frac{1}{2}D^{HUU}$$

$$\frac{1}{2}D^{CL} + \frac{1}{2}D^{KK}$$

$$\frac{1}{2}D^{KK} + \frac{1}{2}D^{HUU}$$

$$\frac{1}{3}D^{CL} + \frac{1}{3}D^{HUU} + \frac{1}{3}D^{KK} \quad (12)$$

We implement the Constrained Dogleg method in the Matlab code CoDoSol using Elliptical trust-region with initial radius of 1. The limiting number of the iterations is set to be 300 and the limiting number of F-evaluations is set to be 1000. The experiment was on 15 problems with dimension between n=2 and n=1000 specified in Table 1.

Different types of constrained systems including systems with solutions both within the feasible region and on the boundary, systems with only lower (upper) bounds and systems with variable components bounded from above and below can be found in this table.

Nonlinear constrained systems come from [13, 15] (problems Pb1 to Pb6), [40] (problems Pb10 to Pb13), chemical equilibrium system given in [17, 32, 33, 41] (Pb7), and nonlinear complementarity problems (NCPs) given in [25, 29] (Pb8, Pb9). While Pb15 [20] comes from nonlinear BVPs. Dealing with large dimensional problems is also critical, these problems need more CPU and, in some cases, they need to be reformed before feeding them into the solver [5]. Problems Pb14, Pb15 are examples of this kind of problems [30]. The starting points for the problems with finite lower and upper bounds are selected by a uniform distribution between l and u i.e. $x_0 = l + 0.25v(u - l), v = 1,2,3$ and for the problems with infinite lower and upper bounds $x_0 = 10^v(1,...,1)^T$ and $x_0 = -10^v(1,...,1)^T, v = 0,1,2, v = 1,2,3$ respectively.

In CoDoSol the trust region size is updated as in [10, 11], the failure criteria and parameter selection is same as [24] where the authors introduced an innovative and efficient method for parameter selection. We tested the algorithm with the scaling matrices and reported the numerical results. Thus, the CoDoSol is tested with the following scaling matrices:

- $D^{CL}$
- $D^{KK}$ with $g = 1$ as suggested in [27]
- $D^{HUU}$ with $p = 1$
- $D^{CON}$ defined by (11)

and $G(x) = D(x)^{-\frac{1}{2}}$. The efficiency of the scaling matrices has been measured by *It*, the number of the iterations and *Fe*, the number of the function evaluations to get convergence (Table 2,3,4).

Bellavia et al showed that CL is superior to KK and KK is superior to HUU on their studied test problems, in our test problems KK is slightly better while HUU shows the poorest performance.

**Table 2.** CoDoSol with different scaling matrices: Number of iterations for the first starting point

| Pb # | KK | | CL | | HUU | | $\frac{1}{3}$KK+$\frac{1}{3}$CL+$\frac{1}{3}$HUU | | $\frac{1}{2}$KK+$\frac{1}{2}$CL | | $\frac{1}{2}$CL+$\frac{1}{2}$HUU | | $\frac{1}{2}$KK+$\frac{1}{2}$HUU | |
|---|---|---|---|---|---|---|---|---|---|---|---|---|---|---|
| | It | Fe | It | Fe | It | Fe | It | Fe | It | Fe | It | Fe | It | Fe |
| 1 | 20 | 27 | 21 | 30 | 36 | 42 | 19 | 24 | 21 | 30 | 21 | 30 | 21 | 30 |
| 3 | 6 | 7 | 6 | 7 | 10 | 12 | 8 | 9 | 6 | 7 | 6 | 7 | 6 | 7 |
| 4 | 5 | 6 | 6 | 7 | 5 | 6 | 5 | 6 | 6 | 7 | 7 | 8 | 7 | 8 |
| 5 | * | | * | | * | | * | | * | | * | | * | |
| 6 | 3 | 4 | 3 | 4 | 3 | 4 | 3 | 4 | 3 | 4 | 3 | 4 | 3 | 4 |
| 7 | 14 | 15 | 13 | 14 | 15 | 16 | 14 | 15 | 13 | 14 | 13 | 14 | 13 | 14 |
| 8 | 25 | 35 | 27 | 38 | * | | 116 | 170 | 28 | 39 | 107 | 145 | 135 | 190 |
| 9 | 7 | 9 | 7 | 9 | 7 | 9 | 7 | 9 | 7 | 9 | 7 | 9 | 7 | 9 |
| 10 | 7 | 10 | 7 | 9 | 7 | 10 | 6 | 9 | 8 | 10 | 7 | 9 | 8 | 10 |
| 11 | 10 | 11 | 10 | 11 | 10 | 11 | 9 | 10 | 9 | 10 | 10 | 11 | 10 | 11 |
| 12 | 6 | 7 | 5 | 6 | 8 | 11 | 6 | 8 | 5 | 6 | 5 | 6 | 5 | 6 |
| 13 | 13 | 14 | 13 | 14 | 22 | 24 | 16 | 18 | 15 | 17 | 12 | 13 | 14 | 16 |
| 14 | 21 | 24 | 21 | 24 | 40 | 42 | 25 | 29 | 21 | 24 | 21 | 24 | 21 | 24 |
| 15 | 8 | 9 | 9 | 11 | 8 | 9 | 9 | 11 | 9 | 11 | 8 | 10 | 8 | 10 |

**Table 3.** CoDoSol with different scaling matrices: Number of iterations for the second starting point

| Pb # | KK | | CL | | HUU | | $\frac{1}{3}$KK+$\frac{1}{3}$CL+$\frac{1}{3}$HUU | | $\frac{1}{2}$KK+$\frac{1}{2}$CL | | $\frac{1}{2}$CL+$\frac{1}{2}$HUU | | $\frac{1}{2}$KK+$\frac{1}{2}$HUU | |
|---|---|---|---|---|---|---|---|---|---|---|---|---|---|---|
| | It | Fe | It | Fe | It | Fe | It | Fe | It | Fe | It | Fe | It | Fe |
| 1 | 6 | 7 | 6 | 7 | 36 | 45 | 19 | 26 | 6 | 7 | 6 | 7 | 6 | 7 |
| 2 | 5 | 6 | 5 | 6 | 6 | 7 | 5 | 6 | 5 | 6 | 5 | 6 | 5 | 6 |
| 4 | 6 | 7 | 6 | 7 | 6 | 7 | 6 | 7 | 6 | 7 | 6 | 7 | 6 | 7 |
| 5 | * | | * | | * | | * | | * | | * | | * | |
| 6 | 5 | 6 | 5 | 6 | 6 | 7 | 5 | 6 | 5 | 6 | 5 | 6 | 5 | 6 |
| 7 | 20 | 21 | 20 | 21 | 31 | 41 | 22 | 21 | 20 | 21 | 20 | 21 | 20 | 21 |
| 8 | 23 | 26 | 23 | 26 | * | | 60 | 68 | 23 | 26 | * | | * | |
| 9 | 9 | 11 | 9 | 11 | 9 | 11 | 9 | 11 | 8 | 10 | 8 | 10 | 9 | 11 |
| 10 | 4 | 5 | 4 | 5 | 4 | 5 | 4 | 5 | 4 | 5 | 4 | 5 | 4 | 5 |
| 11 | 4 | 5 | 4 | 5 | 4 | 5 | 4 | 5 | 4 | 5 | 4 | 5 | 4 | 5 |
| 12 | 1 | 2 | 1 | 2 | 4 | 5 | 1 | 2 | 1 | 2 | 1 | 2 | 1 | 2 |
| 13 | 1 | 2 | 1 | 2 | 4 | 5 | 1 | 2 | 1 | 2 | 1 | 2 | 1 | 2 |
| 14 | 10 | 15 | 10 | 15 | 12 | 13 | 10 | 15 | 10 | 15 | 10 | 15 | 10 | 15 |
| 15 | 6 | 7 | 6 | 7 | 7 | 8 | 6 | 7 | 6 | 7 | 6 | 7 | 6 | 7 |

**Table 4.** CoDoSol with different scaling matrices: Number of iterations for the third starting point

| Pb # | KK | | CL | | HUU | | $\frac{1}{3}$KK+$\frac{1}{3}$CL+$\frac{1}{3}$HUU | | $\frac{1}{2}$KK+$\frac{1}{2}$CL | | $\frac{1}{2}$CL+$\frac{1}{2}$HUU | | $\frac{1}{2}$KK+$\frac{1}{2}$HUU | |
|---|---|---|---|---|---|---|---|---|---|---|---|---|---|---|
| | It | Fe | It | Fe | It | Fe | It | Fe | It | Fe | It | Fe | It | Fe |
| 1 | * | | * | | Err5 | | * | | * | | * | | * | |
| 2 | 4 | 5 | 4 | 5 | 4 | 5 | 4 | 5 | 4 | 5 | 4 | 5 | 4 | 5 |
| 4 | 5 | 6 | 5 | 6 | 5 | 6 | 5 | 6 | 5 | 6 | 5 | 6 | 5 | 6 |
| 5 | 10 | 11 | 10 | 11 | 17 | 18 | 10 | 11 | 10 | 11 | 10 | 11 | 10 | 11 |
| 6 | 7 | 8 | 7 | 8 | 10 | 11 | 7 | 8 | 7 | 8 | 7 | 8 | 7 | 8 |
| 7 | 20 | 21 | 26 | 27 | 27 | 28 | 24 | 25 | 29 | 30 | 27 | 28 | 28 | 29 |
| 8 | 110 | 112 | 113 | 116 | * | | * | | 112 | 115 | * | | * | |
| 9 | * | | * | | * | | * | | * | | * | | * | |
| 10 | 5 | 7 | 5 | 7 | 5 | 6 | 5 | 7 | 5 | 7 | 4 | 6 | 4 | 6 |
| 11 | 7 | 8 | 8 | 9 | 8 | 9 | 8 | 9 | 8 | 9 | 8 | 9 | 8 | 9 |
| 12 | 5 | 6 | 5 | 6 | 6 | 7 | 5 | 6 | 5 | 6 | 5 | 6 | 5 | 6 |
| 13 | 56 | 57 | 55 | 56 | * | | * | | 55 | 56 | 55 | 56 | 55 | 56 |
| 14 | 23 | 26 | 23 | 26 | 42 | 47 | 23 | 26 | 23 | 26 | 23 | 26 | 23 | 26 |
| 15 | 8 | 9 | 7 | 8 | 7 | 8 | 7 | 8 | 7 | 8 | 7 | 8 | 7 | 8 |

The combination of the scaling matrices shows interesting performances. For example, in 80% of the cases $\frac{1}{2}CL + \frac{1}{2}HUU$ shows as good performance as or better performance than the individual matrices.

To compare these seven scaling matrices, we use the performance profiles and to have a more reliable comparison we used Nested Performance Profiles [22], which removed a negative side effect of the performance profiles. In Fig 1 and Fig 3 the computational effort is measured in terms of mean *It* (the mean number of iterations for the three starting points) and mean *Fe* (mean F-evaluations for the three starting points) respectively. In these figures we eliminate the convex coefficients and show $\frac{1}{2}C_1 + \frac{1}{2}C_2$ with $C_1 + C_2$. The black doted lines correspond to the individual matrices and the red line corresponds to the convex combination of the scaling matrices.

By looking at the performance profiles corresponding to CL+HUU we can see that it is efficient in solving about the 75% of the tests and solves about 95% of the tests within a factor 1 from the best solver. CL+HUU is the best scaling matrix for the studied problems.

This can be verified by looking at the figure 2 and figure 4. Figure 2 is based on *It* while Figure 4 is based on *Fe*.

**Table 5.** Projected Affine-Scaling Interior-Point Newton Method with different scaling matrices, Rosenbrock function

| Rosenbrock | KK | CL | HUU | $\frac{1}{3}$KK+$\frac{1}{3}$CL+$\frac{1}{3}$HUU | $\frac{1}{2}$KK+$\frac{1}{2}$CL | $\frac{1}{2}$CL+$\frac{1}{2}$HUU | $\frac{1}{2}$KK+$\frac{1}{2}$HUU |
|---|---|---|---|---|---|---|---|
| Iteration | 3 | 34 | 4 | 4 | 5 | 6 | 3 |
| $\|\|x^k - x^*\|\|$ | 4.0e-18 | 2.45e-13 | 1.53e-18 | 4.81e-17 | 3.85 e-18 | 2.48e-16 | 1.67e-18 |

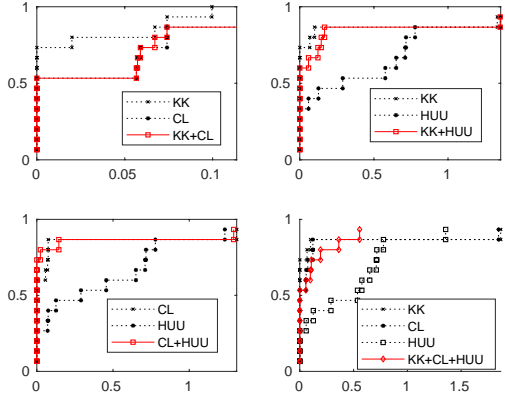

**Figure 1**. Performance Profiles for the number of iterations, basic comparisons

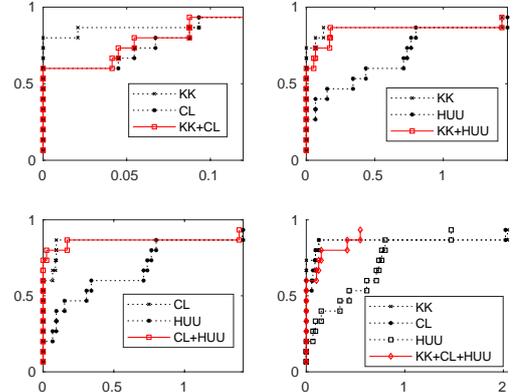

**Figure 3**. Performance Profiles for the F-evaluations, basic comparisons

### B. Experiment 2

**Table 6.** Projected Affine-Scaling Interior-Point Newton Method with different scaling matrices, Wood function

| Wood | KK | CL | HUU | $\frac{1}{3}$KK+$\frac{1}{3}$CL+$\frac{1}{3}$HUU | $\frac{1}{2}$KK+$\frac{1}{2}$CL | $\frac{1}{2}$CL+$\frac{1}{2}$HUU | $\frac{1}{2}$KK+$\frac{1}{2}$HUU |
|---|---|---|---|---|---|---|---|
| Iteration | 3 | 37 | 9 | 7 | 6 | 10 | 5 |
| $\|\|x^k - x^*\|\|$ | 1.0e-18 | 2.95e-14 | 2.57e-18 | 6.34e-17 | 3.41e-18 | 3.23e-16 | 3.12e-18 |

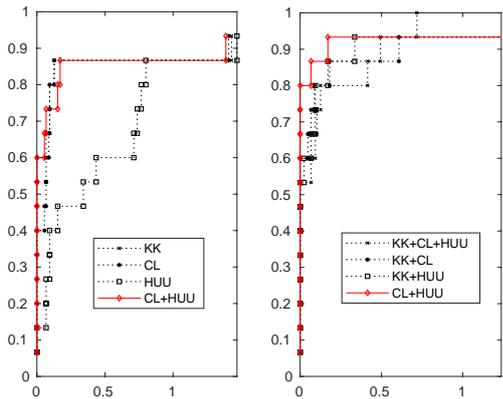

**Figure 2**. Performance Profiles for the number of iterations, the performance of CL+HUU

In this section, we illustrate the local behavior of the different scaling strategies using two standard test problems. We implemented Algorithm "Projected Affine-Scaling Interior-Point Newton Method" [29]. The first test example is the famous Rosenbrock-function:

$$f(x) = 100(x_2 - x_1^2)^2 + (1 - x_1)^2.$$

This function has a unique global minimum at $x^* = (1,1)$. The lower and upper bounds are $l = (0,0)$ and $u = (1,1)$.

To study the local convergence properties, the standard starting point is set to be $x_0 = (0.999, 0.999)$. Table 5 contains the corresponding numerical results for the scaling matrices. The slow linear rate of convergence for CL and significantly better for HUU and KK have been indicated. While for the convex combination we have a fast convergence. This naturally removes the weak point of scaling matrix CL since the number of iterations reduces from 34 to 4.

The second test problem is wood function:

$$f(x) = 100(x_2 - x_1^2) \\ +(1-x_1)^2 + 90(x_4 - x_3)^2 \\ +(1-x_3)^2 + 10(x_2 + x_4 - 2)^2 \\ +0.1(x_2 - x_4)^2.$$

This function admits an unconstrained minimum in $x^* = (1,1,1,1)$. We use the bounds $l = (1,1,1,0.99)$ and $u = (3,3,3,3)$. The starting point is $x^0 = 1.001(1,1,1,1)$. The corresponding numerical results are given in Table 6.

Again, we see that the convergence of the Coleman-Li matrix is rather slow. While according to Bellavia et all it is the fastest one for the medium scale problems. Clearly the convex combination of the CL and KK helps it to mitigate the violation of the strict complementarity assumption that slows down the convergence rate of the affine-scaling Newton method using the Coleman-Li scaling.

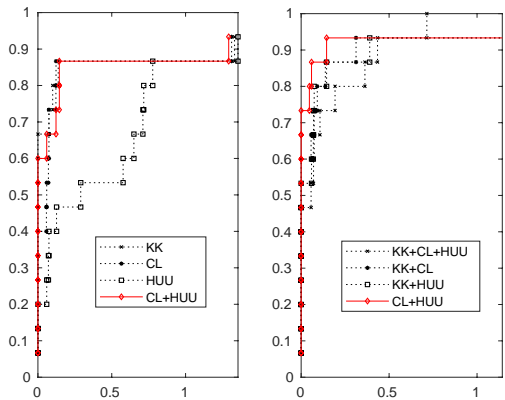

**Figure 4**. Performance Profiles for the F-evaluations, the performance of CL+HUU

## V. CONCLUSION

During our long implementations it has been proved that the performance of the scaling matrices depends on the test problems. Each matrix has its own advantages and disadvantages. For a new problem, there is no way to select the best possible scaling matrix and it has to be done by trial and fail. By changing in the dimension of the problem or changing the parameter values of a problem the previous matrix will not necessarily work as the best one. Also the performance of the scaling matrices highly depends on the algorithm. A matrix demonstrates fast convergence for a specific algorithm and slow convergence for other algorithm. In order to overcome this problem, one can use the convex combination of the scaling matrices. This new class of scaling matrices has the advantages of the individual matrices and in some cases works as the best option.